\documentclass[a4paper,10pt]{article}
\usepackage{amsmath,amssymb}
\usepackage[latin1]{inputenc}
\usepackage[english,frenchb]{babel}

\newcommand{\Dend}{\operatorname{Dend}}
\newcommand{\Com}{\operatorname{Com}}
\newcommand{\Dias}{\operatorname{Dias}}
\newcommand{\Lie}{\operatorname{Lie}}
\newcommand{\Brace}{\operatorname{Brace}}

\newcommand{\QQ}{\mathbb{Q}}
\newcommand{\sym}{\mathfrak{S}}
\newcommand{\tc}{\mathbf{c}}
\newcommand{\ta}{\mathbf{a}}
\newcommand{\tb}{\mathbf{b}}
\newcommand{\fa}{F_\mathbf{a}}
\newcommand{\fb}{F_\mathbf{b}}

\newcommand{\Ind}{\operatorname{Ind}}

\renewcommand{\mod}{\,\operatorname{mod}}

\newtheorem{theorem}{Theorem}[section] 
\newtheorem{proposition}[theorem]{Proposition} 
\newtheorem{conjecture}[theorem]{Conjecture}

\newenvironment{proof}{\begin{trivlist}\item{\bf{Preuve.}}}
  {\hfill\rule{2mm}{2mm}\end{trivlist}}

\title{Le module dendriforme sur le groupe cyclique}
\author{F. Chapoton}
\date{\today}

\begin{document}

\maketitle

\begin{abstract}
  La structure d'opérade anticyclique de l'opérade dendriforme
  donne en particulier une matrice d'ordre $n$ agissant sur l'espace
  engendré par les arbres binaires plans à $n$ feuilles. On
  calcule le polynôme caractéristique de cette matrice. On
  propose aussi une conjecture compatible pour le polynôme
  caractéristique de la transformation de Coxeter du poset de
  Tamari, qui est essentiellement une racine carrée de cette matrice.
\end{abstract}

\selectlanguage{english}

\begin{abstract}
  The structure of anticyclic operad on the Dendriform operad defines
  in particular a matrix of finite order acting on the vector space
  spanned by planar binary trees. We compute its characteristic
  polynomial and propose a (compatible) conjecture for the
  characteristic polynomial of the Coxeter transformation for the
  Tamari lattice, which is mostly a square root of this matrix.
\end{abstract}

\selectlanguage{french}

% MSC 2000 : 18D50 , 05E05, 06A07

\setcounter{section}{-1}

\section{Introduction}

Les arbres binaires plans sont des objets combinatoires très
classiques. Ils ont récemment fait l'objet de recherches
remarquables en combinatoire algébrique, suite à leur apparition
dans les travaux de Loday \cite{livre}. Le point de départ est la
description par Loday d'une opérade basée sur les arbres
binaires plans, nommée opérade dendriforme. A cette opérade
correspond une notion d'algèbre dendriforme. Loday et Ronco ont
ensuite montré \cite{lodayronco} que l'algèbre dendriforme libre
sur un générateur avait une structure d'algèbre de Hopf. De
nombreux travaux
\cite{atree,lodayronco2,aguiarsottile,hnt_cras1,hnt_cras2,hnt} ont
suivi sur différents aspects de ces objets.

L'opérade Dendriforme a en fait une structure plus riche : c'est
une opérade anticyclique \cite{anticyclic}. Ceci donne en
particulier une action du groupe cyclique d'ordre $n+1$ sur l'espace
vectoriel engendré par les arbres binaires plans à $n$ sommets
internes. Cette action reste relativement mystérieuse depuis son
introduction, malgré quelques progrès effectués depuis dans
sa compréhension. On a montré dans \cite{canabull} que l'action
du générateur $\tau$ du groupe cyclique était liée au
carré de la transformation de Coxeter $\theta$ pour le treillis de
Tamari. On a aussi montré dans \cite{moules} que l'action du
générateur $\tau$ sur les arbres binaires plans admettait une
description très simple par le biais d'un plongement des arbres
binaires plans dans les fractions rationnelles.

On calcule ici le polynôme caractéristique du générateur
$\tau$ de cette action cyclique. On propose ensuite une conjecture
pour le polynôme caractéristique de la transformation de Coxeter
$\theta$. On montre en particulier que cette conjecture est compatible
avec le résultat précédent.

Ce travail a son origine dans une correspondance avec Jean Ecalle, qui
a formulé le premier une conjecture sur le polynôme
caractéristique de $\tau$.

Merci à Cédric Bonnafé pour la preuve de la Proposition \ref{injection}.

\section{Rappels et notations}

On note $\Dend$ l'opérade dendriforme. L'espace vectoriel $\Dend(n)$ a
pour base les arbres binaires plans à $n$ sommets internes et pour
dimension le nombre de Catalan  $\tc_n$  défini par
\begin{equation}
  \tc_n = \frac{1}{n+1}\binom{2n}{n}.
\end{equation}

On trouve dans \cite{anticyclic} la description en termes d'opérade
d'une action naturelle du groupe cyclique d'ordre $n$ sur
$\Dend(n-1)$. On a donc un module induit sur le groupe symétrique
$\sym_{n}$.

\begin{proposition}[\cite{anticyclic}, Th. 6.2]
  Le caractère de ce module induit est donné par la fonction
  symétrique
  \begin{equation}
    \label{carac_sym}
    2 \tc_{n-1} p_1^n-\frac{1}{2n} \sum_{d|n} \phi(d) \binom{2n/d}{n/d}p_d^{n/d},
  \end{equation}
  où $\phi$ est l'indicatrice d'Euler et les $p_d$ sont les
  fonctions symétriques ``sommes de puissances''.
\end{proposition}

Cette fonction symétrique décrit donc le module sur le groupe
symétrique $\sym_n$ induit depuis un module $\Dend(n-1)$ sur le groupe
cyclique $C_n$. On cherche à comprendre ce module sur le groupe
cyclique, dont on sait a priori qu'il est défini sur les rationnels et
même sur les entiers.

On voit que la formule (\ref{carac_sym}) se décompose en une somme de
deux termes. Le premier terme $2 \tc_{n-1} p_1^n$ est une somme de $2
\tc_{n-1}$ représentations régulières de $\sym_n$, donc est isomorphe
à l'induite de la somme de $2 \tc_{n-1}$ représentations régulières de
$C_n$. On se concentre donc par la suite sur le deuxième terme.

\section{Induction du groupe cyclique au groupe symétrique}

Soit $C_n$ le groupe cyclique à $n$ éléments. On note $t$ le
générateur. L'algèbre de groupe est $\QQ[t]/(t^n-1)$.

Soit $d$ un diviseur de $n$. L'espace $\QQ[t]/(t^d-1)$ est un module
pour $C_n$, noté $M_{n,d}$. 

Remarque : $M_{n,n}$ est la représentation régulière.

Le caractère de $M_{n,d}$ est le suivant :
\begin{equation}
  \label{chi_m_nd}
  \chi_{n,d}(t^k)=
  \begin{cases}
    d  \text{ si } d|k,\\
    0 \text{ sinon.}
  \end{cases}
\end{equation}

Soit $K_0^{\QQ}(C_n)$ le groupe de Grothendieck de la catégorie des
$\QQ C_n$-modules de type fini. Le groupe $K_0^\QQ(C_n)$ a pour rang
le nombre de diviseurs de $n$ et les classes dans $K_0^{\QQ}(C_n)$ des
modules $M_{n,d}$ pour $d|n$ forment une base de $K_0^{\QQ}(C_n)$
\cite[\S 39, Ex. 2]{curtis_reiner}.

Soit $\sym_n$ le groupe symétrique sur $\{1,\dots,n\}$.

On a une inclusion de groupe de $C_n$ dans $\sym_n$ qui envoie $t$ sur
le grand cycle $c=(1,2,\dots,n)$.

On considère le module induit $M'_{n,d}=\Ind_{C_n}^{\sym_n} M_{n,d}$.

\begin{proposition}
  Le caractère de $M'_{n,d}$ est la fonction symétrique
  \begin{equation}
    \label{chi'}
    \chi'_{n,d}=\frac{d}{n} \sum_{\ell | n/d} \phi(\ell) p_\ell^{n/\ell}.
  \end{equation}
\end{proposition}

\begin{proof}
  Par la formule d'induction des caractères, le caractère
  $\chi'_{n,d}$ est donné (en tant que fonction centrale sur le groupe
  symétrique) par
  \begin{equation}
    \chi'_{n,d}(\sigma)=\frac{1}{n}\sum_{{\tau \in \sym_n}\atop
      {\tau^{-1}\sigma\tau\in C_n}} \chi_{n,d}(\tau^{-1}\sigma\tau),
  \end{equation}
  pour toute permutation $\sigma$. En utilisant (\ref{chi_m_nd}), ceci
  se simplifie en
  \begin{equation}
    \frac{d}{n}\sum_{d|k}\sum_{{\tau \in \sym_n}\atop
     { \tau^{-1}\sigma\tau=c^k}} 1.
  \end{equation}
  On note que le type cyclique de $c^k$ est $(\frac{n}{n \wedge k})^{n
    \wedge k}$. On traduit alors le caractère $\chi'_{n,d}$ en termes de
  fonctions symétriques :
  \begin{equation}
    \chi'_{n,d}=\frac{d}{n} \sum_{{k=1\dots n} \atop {d|k}} (p_{\frac{n}{n
        \wedge k}})
    ^{n \wedge k}.
  \end{equation}
  On remplace ensuite $k$ par $d k$ et $n$ par $d m$ :
  \begin{equation}
    \frac{d}{n} \sum_{k=1}^{m} (p_{\frac{m}{m \wedge k}})^{(m \wedge
      k)d}.
  \end{equation}
  On trouve la formule attendue en regroupant les termes.
\end{proof}

\begin{proposition}
  \label{injection}
  L'induction $\Ind_{C_n}^{\sym_n}$ de $K_0^\QQ(C_n)$ dans
  $K_0^\QQ(\sym_n)$ est une application linéaire injective.
\end{proposition}

\begin{proof}
  Comme les classes des modules $M_{n,d}$ forment une base de
  $K_0^\QQ(C_n)$, il suffit de montrer que les classes des modules
  induits $M'_{n,d}$ dans $K_0^\QQ(\sym_n)$ sont linéairement
  indépendantes. Ceci résulte immédiatement de la formule
  (\ref{chi'}), par triangularité par rapport à l'ordre partiel
  défini sur l'ensemble des diviseurs de $n$ par la
  divisibilité.
\end{proof}

\section{Description du module dendriforme}

On définit une suite d'entiers positifs $\ta_n$ pour $n \geq 1$ par la
formule suivante :
\begin{equation}
  \ta_n=\frac{1}{2n}\sum_{d|n} \mu(n/d) \binom{2d}{d},
\end{equation}
où $\mu$ est la fonction de Möbius.

A priori, les $\ta_n$ définis ainsi sont des rationnels. Il faut
justifier que ce sont bien des entiers positifs. Pour cela, on va
faire un calcul plus fin avec des fonctions symétriques. On note
$\circ$ le pléthysme des fonctions symétriques.

Considérons la fonction symétrique $\Lie$ :
\begin{equation}
  \Lie=\sum_{n \geq 1} \frac{1}{n}\sum_{d|n} \mu(d) p_d ^{n/d}.
\end{equation}
C'est la fonction symétrique associée à l'opérade $\Lie$, donc c'est
une somme positive de fonctions de Schur.

Considérons la fonction symétrique $\Brace$ :
\begin{equation}
  \Brace=\sum_{n \geq 1} \frac{1}{n}\binom{2n-2}{n-1} p_1
  ^{n}=\frac{1-\sqrt{1-4 p_1}}{2}.
\end{equation}
C'est aussi clairement une somme positive de fonctions de Schur.

\begin{proposition}
  \label{lie_brace}
  Le pléthysme $\Lie \circ \Brace$ est la fonction symétrique
  \begin{equation}
    \label{liebrace}
    \sum_{n \geq 1} \frac{1}{2n}\sum_{d|n} \mu(d) \binom{2n/d}{n/d} p_d ^{n/d}.
  \end{equation}
\end{proposition}
\begin{proof}
  Par définition du pléthysme, il faut calculer
  \begin{equation}
    \sum_{n \geq 1} \frac{1}{n}\sum_{d|n} \mu(d)
    \left( \frac{1-\sqrt{1-4 p_d}}{2}\right)^{n/d}.
  \end{equation}
  En échangeant les sommations, ceci devient
  \begin{equation}
    - \sum_{d\geq 1} \frac{\mu(d)}{d} \log\left(\frac{1+\sqrt{1-4 p_d}}{2}\right).
  \end{equation}
  En utilisant alors le développement de Taylor (\ref{eq_A3}) (Voir
  Appendice) :
  \begin{equation}
    -\log\left(\frac{1+\sqrt{1-4 u}}{2}\right)
=\sum_{n \geq 1} \frac{1}{2n}\binom{2n}{n} u^n,
  \end{equation}
  on obtient donc
  \begin{equation}
    \sum_{d \geq 1}  \frac{\mu(d)}{d} \sum_{k \geq 1} \frac{1}{2k}\binom{2k}{k} p_d^k,
  \end{equation}
  qui se réduit facilement au résultat voulu.
\end{proof}

Comme composée pour le pléthysme de sommes positives de
fonctions de Schur, la fonction $\Lie \circ \Brace$ est aussi une
somme positive de fonctions de Schur. La dimension des invariants dans
$\Lie \circ \Brace$, obtenue en posant $p_d=1$ pour tout $d$ dans les
composantes homogènes de la formule (\ref{liebrace}), est
exactement la suite $\ta_n$ ; ces nombres sont donc bien des entiers
positifs.

Voici les premiers termes de la suite $\ta_n$ pour $n \geq 1$ :
\begin{equation}
   1,1,3,8,25,75,245,800,2700,9225,\dots
\end{equation}

La suite $\ta_n$ a une propriété remarquable.
\begin{proposition}
  On a
  \begin{equation}
    \label{prod_ta}
    \fa(x)=\prod_{n \geq 1}(1-x^n)^{-\ta_n}= \frac{1-\sqrt{1-4x}}{2x}.
  \end{equation}
\end{proposition}

\begin{proof}
  Soit $\Com$ la fonction symétrique associée à l'opérade $\Com$ et
  correspondant à la somme des modules triviaux sur les groupes
  symétriques. Elle vérifie $(1+\Com) \circ \Lie= \frac{1}{1-p_1}$.
  Par conséquent, on a
  \begin{equation}
    (1+\Com) \circ (\Lie \circ \Brace) 
= \left( \frac{1}{1-p_1} \right) \circ \Brace.
  \end{equation}
  En regardant les invariants de part et d'autre, on en déduit la proposition.
\end{proof}

On utilise maintenant cette suite d'entiers $\ta_n$ pour décrire le
module dendriforme $\Dend(n-1)$ sur le groupe cyclique $C_n$.

\begin{proposition}
  Pour tout $n \geq 1$, on a l'égalité suivante :
  \begin{equation}
    \sum_{d|n} \ta_d \chi'_{n,d}=\frac{1}{2n} \sum_{d|n} \phi(d) \binom{2n/d}{n/d}p_d^{n/d}.
  \end{equation}
\end{proposition}

\begin{proof}
  On calcule la somme pour $n\geq 1$ du membre de gauche. En utilisant
  l'expression (\ref{chi'}) de $\chi'_{n,d}$, on obtient
  \begin{equation}
    \sum_{n \geq 1} \sum_{d|n} \frac{1}{2d}\sum_{k|d} \mu(d/k)
    \binom{2k}{k} \frac{d}{n} \sum_{\ell | n/d} \phi(\ell) p_\ell^{n/\ell},
  \end{equation}
  ce qui s'écrit encore
  \begin{equation}
    \sum_{n \geq 1}\sum_{d|n} \sum_{k|d} \sum_{\ell | n/d}  \frac{1}{2n}  \mu(d/k)
    \binom{2k}{k} \phi(\ell) p_\ell^{n/\ell}. 
  \end{equation}
  
  On introduit de nouvelles variables de sommation $i$ et $j$ en
  posant $n=i j k \ell$ et $d=i k$. En remplaçant les sommations sur
  $n$ et $d$ par des sommations sur $i$ et $j$, on obtient
  \begin{equation}
    \sum_{i,j,k,\ell} \frac{1}{2 i j k \ell} \mu(i) \phi(\ell) \binom{2k}{k} p_\ell^{i j k}. 
  \end{equation}
  
  En utilisant le pléthysme des fonctions symétriques, on peut
  factoriser cette expression comme suit :
  \begin{equation}
    \left(\sum_i \mu(i) p_i/i \right) \circ \left(\sum_j p_j/j\right) \circ \left(\sum_{k,\ell}
      \frac{1}{2 k \ell} \phi(\ell) \binom{2k}{k} p_\ell^{ k}\right).
  \end{equation}
  
  Comme les deux premiers termes sont dans un sous-groupe des
  fonctions symétriques pour le pléthysme qui est isomorphe au groupe
  commutatif des séries de Dirichlet pour le produit, ces termes sont
  inverses l'un de l'autre. On obtient donc
  \begin{equation}
    \sum_{k,\ell}
    \frac{1}{2 k \ell} \phi(\ell) \binom{2k}{k} p_\ell^{ k}.
  \end{equation}
  Il est facile de voir que le terme de degré $n$ de cette somme est
  exactement le résultat attendu.
\end{proof}

On a donc montré que l'induite du module virtuel
\begin{equation}
  2\tc_{n-1} M_{n,n} - \bigoplus_{d|n} \ta_d M_{n,d}
\end{equation}
a le même caractère que l'induite du module $\Dend(n-1)$.

On sait aussi que le module $\Dend(n-1)$ est défini sur les rationnels
et même sur les entiers. On sait par ailleurs que les modules
$M_{n,d}$ sont aussi définis sur les rationnels.

Par la proposition \ref{injection}, on en déduit 

\begin{theorem}
  \label{theo_principal}
  Le module $\Dend(n-1)$ a pour caractère 
  \begin{equation}
    2\tc_{n-1} \chi_{n,n} - \bigoplus_{d|n} \ta_d \chi_{n,d},
  \end{equation}
  et le polynôme caractéristique du générateur $t$ de $C_n$ est
  \begin{equation}
    \label{poly_tau}
    \frac{(t^{n}-1)^{2 \tc_{n-1}}}{\prod_{d| n}(t^d-1)^{\ta_d}}.
  \end{equation}
\end{theorem}

Ce théorème incite à penser que le module $\Dend(n-1)$ doit admettre
une résolution de la forme
\begin{equation}
0 \longrightarrow \bigoplus_{d|n} \ta_d M_{n,d} \longrightarrow 2\tc_{n-1} M_{n,n}
  \longrightarrow \Dend(n-1) \longrightarrow 0.
\end{equation}

On a une application évidente de $c_{n-1} M_{n,n}$ dans
$\Dend(n-1)$, induite par l'identité de $\Dend(n-1)$. Pour définir un
morphisme de $2 c_{n-1} M_{n,n}$ dans $\Dend(n-1)$, il faudrait une
autre application de $c_{n-1} M_{n,n}$ dans $\Dend(n-1)$. On peut
supposer qu'elle doit provenir d'une involution sur $\Dend(n-1)$.

\bigskip

Remarque : dans le cas similaire mais plus simple de l'opérade
$\Dias$, dont le caractère anticyclique est donné par
$M'_{n,n}-M'_{n,1}$, on retrouve la description connue du module
$\Dias(n-1)$ comme quotient de $M_{n,n}$ par $M_{n,1}$. On a une suite
exacte courte
\begin{equation}
0 \longrightarrow M_{n,1} \longrightarrow M_{n,n}
  \longrightarrow \Dias(n-1) \longrightarrow 0.
\end{equation}

\section{Transformation de Coxeter}

Le treillis de Tamari \cite{tamari} est un ordre partiel sur les
arbres binaires plans à $n$ feuilles. Soit $C$ la matrice de cet
ordre :
\begin{equation}
  C_{x,y}=
  \begin{cases}
    1 \text{ si } x\leq y,\\
    0 \text{ sinon}.
  \end{cases}
\end{equation}
La matrice $\theta=-C({}^{t}C^{-1})$ est appelée la transformation
de Coxeter du poset de Tamari.

On a montré dans \cite{canabull} la relation suivante.
\begin{proposition}[\cite{canabull}, Th. 6.1]
  Si $\tau$ est la matrice d'ordre $n$ considérée précédemment et
  $\theta$ la transformation de Coxeter du treillis de Tamari, alors
  on a
  \begin{equation}
    \tau=(-1)^{n+1} \theta^2.
  \end{equation}
\end{proposition}

Il est donc naturel de se demander si il existe une description simple
du polynôme caractéristique de $\theta$. On propose ci-dessous une
conjecture pour ce polynôme. Il est nécessaire de distinguer les cas
$n$ pair et $n$ impair.

On commence par introduire une suite d'entiers relatifs $\lambda(n)$.
\begin{equation}
   \lambda(n)=(-1)^{\binom{n}{2}} \binom{n-1}{\lceil (n-1)/2 \rceil}. 
\end{equation}

On définit ensuite une suite d'entiers relatifs $\tb_{n}$ par
inversion de Möbius :
\begin{equation}
  \tb_{n}=\frac{1}{n}\sum_{d|n} \mu(d)\lambda(n/d).
\end{equation}

A priori, les $\tb_n$ sont des rationnels. Il faut justifier que ce
sont bien des entiers relatifs. On procède comme pour la suite
$\ta_n$, en utilisant des fonctions symétriques.

On introduit la fonction symétrique
\begin{equation}
  Z=p_1+\sum_{n \geq 0} \frac{(-1)^{n+1}}{n+1}\binom{2n}{n}  p_1^{2n+2} 
=\frac{1+2 p_1-\sqrt{1+4 p_1^2}}{2}.
\end{equation}

\begin{proposition}
  La fonction symétrique $\Lie \circ Z$ est
  \begin{equation}
    \sum_{n \geq 1} \frac{1}{n} \sum_{d|n} \mu(d)\lambda(n/d) p_d^{n/d}.
  \end{equation}
\end{proposition}

\begin{proof}
  La preuve est similaire à celle de la Prop. \ref{lie_brace}. Par
  définition du pléthysme, il faut calculer
  \begin{equation}
    \sum_{n \geq 1} \frac{1}{n}\sum_{d|n} \mu(d)
    \left( \frac{1+2 p_d-\sqrt{1+4 p_d^2}}{2}\right)^{n/d}.
  \end{equation}
  En échangeant les sommations, on obtient
  \begin{equation}
    - \sum_{d\geq 1} \frac{\mu(d)}{d} \log\left(\frac{1-2 p_d + \sqrt{1+4 p_d^2}}{2}\right).
  \end{equation}

  On utilise alors le développement de Taylor (\ref{eq_A6}) (Voir
  Appendice) :
  \begin{equation}
    -\log\left(\frac{1-2u+\sqrt{1+4u^2}}{2}\right)
    =\sum_{n \geq 1} \lambda(n) u^n/n.
  \end{equation}

  On obtient donc
  \begin{equation}
    \sum_{d \geq 1}  \frac{\mu(d)}{d} \sum_{k \geq 1} \frac{\lambda(k)}{k} p_d^k,
  \end{equation}
  qui se réduit facilement au résultat voulu.
\end{proof}

Comme $\Lie$ et $Z$ sont des combinaisons linéaires entières de
fonctions de Schur, la fonction symétrique $\Lie \circ Z$ l'est aussi.
Ceci montre que les $\tb_n$ sont bien des entiers relatifs, car ce
sont les coefficients des invariants dans cette fonction symétrique.

Voici les premiers termes de la suite $\tb_n$ pour $n \geq 1$ :
\begin{equation}
   1, -1, -1, 1, 1, -1, -3, 4, 8, -13, -23, 39, 71, -121,\dots
\end{equation}

On remarque que les signes des $\tb_n$ semblent suivre un motif
régulier, le même que pour la suite $\lambda(n)$. En fait, il
semble même que la fonction symétrique $\Lie \circ Z$ soit une
somme de fonctions de Schur dont les signes sont soit tous positifs,
soit tous négatifs selon la valeur de $n$ modulo $4$. On a
vérifié ceci pour $n \leq 12$.

\begin{proposition}
  On a
  \begin{equation}
    \label{prod_tb}
    \fb(x)=\prod_{n \geq 1}(1-x^n)^{-\tb_n}= \frac{-1+2x+\sqrt{1+4 x^2}}{2x}.
  \end{equation}
\end{proposition}
\begin{proof}
  On utilise l'identité
  \begin{equation}
    (1+\Com) \circ (\Lie \circ Z)= \frac{1}{1-p_1} \circ Z.
  \end{equation}
  En prenant les invariants, on trouve l'égalité voulue.
\end{proof}

On utilise maintenant la suite de nombres entiers relatifs $\tb_n$ pour
proposer une description de $\theta$.

\begin{conjecture}
  \label{conj_pair}
  Pour $n$ pair, le polynôme caractéristique de $\theta$ est
  \begin{equation}
    \frac{(t^{2n}-1)^{\tc_{n-1}}}{\prod_{d| 2n}(t^d-1)^{\tb_d}}.
  \end{equation}
\end{conjecture}

Pour donner une conjecture dans le cas $n$ impair, on définit une
autre suite d'entiers relatifs $\tb'_n$ à partir de la suite $\tb_n$ :
\begin{equation}
  \label{fromb2b}
  \tb'_n=
  \begin{cases}
    \tb_n &\text{ si } n = 1 \mod 2,\\
    -\tb_{n}-\tb_{n/2}  &\text{ si } n = 0 \mod 2. 
  \end{cases}  
\end{equation}

Voici les premiers termes de la suite $\tb'_n$ pour $n \geq 1$ :
\begin{equation}
  1, 0, -1, 0, 1, 2, -3, -5, 8, 12, -23, -38, 71, 124, \dots
\end{equation}

\begin{conjecture}
  \label{conj_impair}
  Pour $n$ impair, le polynôme caractéristique de $\theta$ est
  \begin{equation}
    \frac{(t^{2n}-1)^{\tc_{n-1}}}{\prod_{d| 2n}(t^d-1)^{\tb'_d}}.
  \end{equation}
\end{conjecture}

\section{Comparaison entre conjectures et théorème}

On montre ici que les conjectures \ref{conj_pair} et
\ref{conj_impair} sont compatibles avec (et impliquent) le théorème
\ref{theo_principal} décrivant le polynôme caractéristique pour
$\tau$.

Le polynôme caractéristique du carré d'une matrice $M$ d'ordre fini
est obtenu par la substitution suivante dans celui de $M$ :
\begin{equation}
  (t^d-1) \mapsto
  \begin{cases}
    t^d-1 &\text{ si } d=1 \mod 2,\\
    (t^{d/2}-1)^2 &\text{ si } d=0 \mod 2.
  \end{cases}
\end{equation}
Le polynôme caractéristique de l'opposé d'une matrice $M$ d'ordre fini
est obtenu par la substitution suivante dans celui de $M$ :
\begin{equation}
  (t^d-1) \mapsto
  \begin{cases}
    (t^{2d}-1)/(t^d-1) &\text{ si } d=1 \mod 2,\\
    (t^{d}-1) &\text{ si } d=0 \mod 2.
  \end{cases}
\end{equation}

Supposons d'abord $n$ impair et considérons la conjecture
\ref{conj_impair}. Dans ce cas, on a $\tau=\theta^2$. On obtient donc
l'expression suivante pour le polynôme caractéristique de $\tau$ :
\begin{equation}
   {(t^{n}-1)^{2 \tc_{n-1}}}{\prod_{{d| 2n}\atop d=1 \mod 2}(t^d-1)^{-\tb'_d}}{\prod_{{d| 2n}\atop d=0 \mod 2}(t^{d/2}-1)^{-2\tb'_d}}.
\end{equation}
Ceci se ré-écrit
\begin{equation}
   (t^{n}-1)^{2 \tc_{n-1}} {\prod_{d| n}(t^d-1)^{-\tb'_d-2\tb'_{2d}}}.
\end{equation}
Pour identifier ceci à la formule (\ref{poly_tau}), il faut donc avoir
\begin{equation}
  \ta_d=\tb'_d+2\tb'_{2d}
\end{equation}
pour tous les entiers $d$ impairs. En utilisant les relations
(\ref{fromb2b}), ceci est équivalent à
\begin{equation}
  \label{condi_impair}
  \ta_d=-\tb_d-2\tb_{2d}
\end{equation}
pour tous les entiers $d$ impairs.

Supposons maintenant $n$ pair et considérons la conjecture
\ref{conj_pair}. Dans ce cas, on a $\tau=-\theta^2$. On obtient donc
l'expression suivante pour le polynôme caractéristique de $\tau$ :
\begin{equation}
     (t^{n}-1)^{2 \tc_{n-1}}
     \prod_{{d| n}\atop d=0 \mod 2}(t^d-1)^{-2\tb_{2d}}
     \prod_{{d| n}\atop d=1 \mod 2}\frac{(t^{2d}-1)^{-2\tb_{2d}}}{(t^{d}-1)^{-2\tb_{2d}}}
     \prod_{{d| n}\atop d=1 \mod 2}\frac{(t^{2d}-1)^{-\tb_d}}{(t^{d}-1)^{-\tb_d}}.
\end{equation}
Pour identifier cette expression avec la formule (\ref{poly_tau}), on
distingue selon la valeur de $d$ modulo $4$. On obtient les conditions
suivantes :
\begin{equation}
 \label{crux}
  \ta_d=
  \begin{cases}
   -2 \tb_{2d}-\tb_d &\text{ si }d=1 \mod 2,\\
  2 \tb_{2d} &\text{ si }d=0 \mod 4,\\
  2 \tb_{2d}+2\tb_d+\tb_{d/2} &\text{ si }d=2 \mod 4.
  \end{cases}
\end{equation}
On remarque que la condition (\ref{condi_impair}) obtenue plus haut
dans le cas $n$ impair fait aussi partie des trois conditions
ci-dessus.

\begin{proposition}
  On a la relation
  \begin{multline}
    \prod_{n} (1-x^n)^{-\ta_n}=\\
    \prod_{n = 0 \mod 4} (1-x^n)^{-2\tb_{2n}}
    \prod_{n = 2 \mod 4} (1-x^n)^{-2\tb_{2n}-2\tb_n-\tb_{n/2}}
    \prod_{n = 1 \mod 2} (1-x^n)^{2\tb_{2n}+\tb_n}.
  \end{multline}
\end{proposition}

\begin{proof}
Calculons le second membre :
\begin{equation}
   \prod_{n = 0 \mod 2} (1-x^{2n})^{-2\tb_{4n}}
  \prod_{n = 1 \mod 2} (1-x^{2n})^{-2\tb_{4n}-2\tb_{2n}-\tb_{n}}
  \prod_{n = 1 \mod 2} (1-x^n)^{2\tb_{2n}+\tb_n}. 
\end{equation}
On coupe en deux le facteur central et on regroupe :
\begin{equation}
    \prod_{n } (1-x^{2n})^{-2\tb_{4n}}
  \prod_{n = 1 \mod 2} (1+x^{n})^{-2\tb_{2n}-\tb_n}.  
\end{equation}
On ré-écrit le premier facteur et on coupe en deux le second :
\begin{equation}
    \prod_{n = 0 \mod 2} (1-x^{n})^{-2\tb_{2n}}
  \prod_{n = 1 \mod 2} (1+x^{n})^{-2\tb_{2n}}  
  \prod_{n = 1 \mod 2} (1+x^{n})^{-\tb_n}.  
\end{equation}
On regroupe les deux premiers facteurs et on ré-écrit le troisième :
\begin{equation}
    \prod_{n } (1-(-x)^{n})^{-2\tb_{2n}}
  \prod_{n = 1 \mod 2} (1-(-x)^{n})^{-\tb_n}.  
\end{equation}
On complète le second facteur :
\begin{equation}
    \prod_{n } (1-(-x)^{n})^{-2\tb_{2n}}
    \prod_{n } (1-(-x)^{2n})^{\tb_{2n}} 
  \prod_{n } (1-(-x)^{n})^{-\tb_n}.  
\end{equation}
Le troisième facteur est $\fb(-x)$, la fonction $\fb$ étant définie
par (\ref{prod_tb}). On regroupe les deux premiers facteurs et on
simplifie :
\begin{equation}
    \prod_{n } \left(\frac{1-(-x)^{n}}{1+(-x)^{n}}\right)^{-\tb_{2n}} \fb(-x).
\end{equation}
On pose $x=-z^2$ et on choisit les signes avec soin :
\begin{equation}
    \prod_{n } \left(\frac{1-z^{2n}}{1+(-z)^{2n}}\right)^{-\tb_{2n}} \fb(z^2),
\end{equation}
soit encore
\begin{equation}
    \prod_{n =0 \mod 2 }
    \left(\frac{1-z^{n}}{1+(-z)^{n}}\right)^{-\tb_{n}} \fb(z^2).
\end{equation}
On peut alors compléter le produit sans introduire de nouveaux termes :
\begin{equation}
    \prod_{n } \left(\frac{1-z^{n}}{1+(-z)^{n}}\right)^{-\tb_{n}} \fb(z^2),
\end{equation}
soit enfin
\begin{equation}
  \fb(z)
  \prod_{n } (1-(-z)^{2n})^{\tb_{n}}
  \prod_{n } (1-(-z)^{n})^{-\tb_{n}}
  \fb(z^2).
\end{equation}
On obtient donc
\begin{equation}
  \fb(z)(1/\fb(z^2))\fb(-z)\fb(z^2)=\fb(z)\fb(-z),
\end{equation}
ce qui vaut bien
\begin{equation}
  \frac{1-\sqrt{1+4 z^2}}{-2z^2}=    \frac{1-\sqrt{1-4 x}}{2x}.
\end{equation}
Ceci est bien la fonction $\fa(x)$ définie en
(\ref{prod_ta}), comme attendu.
\end{proof}

Remarque : on a montré au passage la relation
\begin{equation}
  \fa(-z^2)=\fb(z)\fb(-z).
\end{equation}

\section*{Appendice}

On rappelle quelques développements de Taylor.

Le premier est classique et facile :
\begin{equation}
  \label{eq_A1}
  \sum_{n \geq 0} \binom{2n}{n} x^n = \frac{1}{\sqrt{1-4x}}.
\end{equation}

On déduit de (\ref{eq_A1}) en utilisant l'opérateur $x \partial_x$ :
\begin{equation}
  \label{eq_A2}
  \sum_{n \geq 1} \frac{1}{n}\binom{2n-2}{n-1} x^n = \frac{1-\sqrt{1-4x}}{2}.
\end{equation}

On déduit aussi de (\ref{eq_A1}) en utilisant l'opérateur $x
\partial_x$ :
\begin{equation}
  \label{eq_A3}
  \sum_{n \geq 1} \frac{1}{2n}\binom{2n}{n} x^n = -\log\left(\frac{1+\sqrt{1-4x}}{2}\right).
\end{equation}

On déduit de (\ref{eq_A2}) en remplaçant $x$ par $-x^2$ :
\begin{equation}
  \label{eq_A4}
  \sum_{n \geq 0} \frac{(-1)^{n+1}}{n+1}\binom{2n}{n} x^{2n+1} = \frac{1-\sqrt{1+4x^2}}{2x}.
\end{equation}

On démontre ensuite
\begin{equation}
  \label{eq_A5}
  \sum_{n \geq 1} \lambda_n x^n = \frac{x}{\sqrt{1+4x^2}}-\frac{1}{2}\left(1-\frac{1}{\sqrt{1+4x^2}}\right),
\end{equation}
en séparant les puissances paires et impaires de $x$ et en utilisant
la relation
\begin{equation}
  \binom{2n}{n} = 2 \binom{2n-1}{n}.
\end{equation}

On déduit ensuite de (\ref{eq_A5}) en utilisant l'opérateur $x \partial_x$ :
\begin{equation}
  \label{eq_A6}  
  \sum_{n \geq 1} \frac{\lambda_n}{n} x^n=-\log\left(\frac{1-2x+\sqrt{1+4x^2}}{2}\right).
\end{equation}

\bibliographystyle{alpha}
\bibliography{poly_carac_dend}

\end{document}